\documentclass[final]{siamart171218}

\usepackage{amssymb}
\usepackage[utf8]{inputenc}
\usepackage{cite}
\usepackage{comment}
\usepackage{mathtools}
\usepackage{algorithmic}
\usepackage{booktabs}
\usepackage{multirow}
\usepackage{siunitx}
\usepackage{capt-of}
\usepackage{mathtools}
\usepackage{tikz}
\usepackage{microtype}
\usepackage{xcolor}

\renewcommand{\Im}{\mathop{\mathrm{Im}}%
}
\renewcommand{\Re}{\mathop{\mathrm{Re}}%
}
\newcommand{\Res}{\mathop{\mathrm{Res}}%
}

\def\resultssetup{
  \sisetup{
    table-number-alignment = center,
    table-figures-integer = 1,
    table-figures-decimal = 2,
    table-figures-exponent = 2,
    exponent-product = \cdot,
    table-sign-exponent
  }
  \setlength{\tabcolsep}{4pt}
}

\newenvironment{resultstable}
{
  \begin{table}
  \centering
  \resultssetup
}
{
  \end{table}
}

\def\ifunc{\mathcal{M}}

\allowdisplaybreaks[1]

\title{Conformal Mapping via a Density Correspondence for
    the Double-Layer Potential}

\author{Matt Wala%
  \thanks{Department of Computer Science, University of Illinois at Urbana-Champaign,
    201 N. Goodwin Ave, Urbana, IL 61801
    (\email{wala1@illinois.edu}, \email{andreask@illinois.edu}).}
  \and
  Andreas Klöckner\footnotemark[1]
}

\def\Emap{f^+}
\def\pEmap{{(\Emap)}}

\def\Imap{f^-}

\def\Dom{\Omega}
\def\Idom{\Dom^-}

\def\Ebcorr{\theta^+}

\def\Ibcorr{\theta^-}

\def\Edom{\Omega^+}

\definecolor{imapcolor}{RGB}{255,127,14}
\definecolor{emapcolor}{RGB}{31,119,180}

\newcommand{\dlk}[2]{\hat{n} \cdot \nabla_{#2} \log {\lvert #2 - #1 \rvert}}

\def\nevalpts{36}

\def\squareorder{9}
\def\pointspersquareinitial{36}
\def\pointspersquarerefinement{72}
\def\pointspersquarefinal{3060}
\def\pointspersnowflakeinitial{1728}
\def\pointspersnowflakerefinement{3456}
\def\pointspersnowflakefinal{136512}
\def\nwhalesourcepoints{2^{13}}

\begin{document}

\maketitle

\begin{abstract}
  We derive a representation formula for harmonic polynomials and Laurent
  polynomials in terms of densities of the double-layer potential on bounded
  piecewise smooth and simply connected domains. From this result, we obtain a
  method for the numerical computation of conformal maps that applies to both
  exterior and interior regions. We present analysis and numerical experiments
  supporting the accuracy and broad applicability of the method.
\end{abstract}



\section{Introduction}

This paper presents an integral equation method for numerical conformal mapping,
using an integral equation based on the \emph{Faber polynomials} (on
the interior) and their counterpart, \emph{Faber-Laurent polynomials}
(on the exterior). Our method is applicable to computing the conformal
map from the interior and exterior of domains bounded by a piecewise
smooth Jordan curve $\Gamma$ onto the interior/exterior of the unit
disk.  Like most techniques for conformal mapping, this method relies
on computing a boundary correspondence function between $\Gamma$ and
the boundary of the target domain. From the boundary correspondence,
the mapping function can be derived via a Cauchy integral~\cite[p.\
381]{henrici}.

The numerical construction of a function that maps the exterior of a simply
connected region conformally onto the exterior of some other region arises in a
number of applications including fluid mechanics~\cite[Ch.\ 4.5]{carrier},
the generation of finite element meshes for problems in fracture mechanics~\cite{tsamasphyros},
the design of optical media~\cite{leonhardt},
the analysis of iterative methods~\cite{starke,driscoll-potential,eiermann}, and
the solution of initial value problems~\cite{moret}. The exterior mapping function's
close relation to the Faber polynomials~\cite{curtiss, suetin} enables a number
of the latter applications. Complex analytic functions defined inside a Jordan
domain admit a near-optimal polynomial expansion in a basis of Faber
polynomials~\cite{curtiss,suetin,ellacott-faber-methods}. By means of the
\emph{Faber transform}, this can be exploited numerically to perform polynomial or
rational approximation~\cite{ellacott-rational-approximation}. A key step of
this approximation procedure is the evaluation of the Faber polynomials
themselves, which, given a numerically derived exterior boundary correspondence,
can be achieved with an FFT-based method~\cite{ellacott} or by
applying Lemma~\ref{lem:happy-lemma} in this paper.

As a technical tool, we introduce a representation formula of complex-valued
harmonic polynomials and Laurent polynomials in terms of the double-layer potential.
The double-layer potential with a complex-valued density $\sigma$ is
given by
\[
\mathcal{D} \sigma(x) =
- \frac{1}{2 \pi} \int_\Gamma \sigma(y) \, \dlk{x}{y} \, ds(y),
\quad x \notin \Gamma.
\]

We study the density functions $\sigma$ that give rise to polynomials on the
\emph{interior} in terms of the \emph{exterior} Riemann map
$\Emap$. Specifically, in this paper, we prove that the images of the $m$th
powers of the values of the exterior map ${(\Emap)}^m$ used as a density under
$\mathcal{D}$ are scaled Faber polynomials of degree $m$. Since the Faber
polynomials form a basis for all complex polynomials, this result characterizes
all densities that give rise to polynomials under the double-layer operator. On
the exterior domain, we study the representation of complex Laurent
polynomials. We find, analogously to the interior case, that the images of the
$m$th powers of values of the interior map $(\Imap)^m$ under $\mathcal{D}$ are
scaled Faber-Laurent polynomials on the exterior domain, where $\Imap$ is the
interior mapping function. In both the polynomial and the Laurent polynomial
case, letting $m = 1$ leads to a uniquely solvable integral equation from which
the boundary correspondence of the interior or exterior map may be recovered by
the solution of an integral equation identical to that of an appropriate Dirichlet problem
of the Laplace equation.

Furthermore, we demonstrate how a Nyström
discretization~\cite{nystroem_praktische_1930} using high-order accurate
quadrature rules achieves high-order accuracy in the computed density. Our
method is of practical interest because of the ready availability of fast
solvers for second kind equations involving the double-layer potential. In the
numerical examples of this paper we use an accelerated solver consisting of
GMRES~\cite{saad_gmres_1986} with the required matrix-vector product driven by
the Fast Multipole Method~\cite{carrier:fmm}. Since the number of GMRES
iterations in our scheme does not depend on the mesh resolution, this solution
scheme has an overall complexity of $O(n \log n)$, where $n$ is the number of
discretization points on the boundary.

Compared with the second-kind integral equation formulations in the existing
literature, our method is perhaps operationally the most similar to
Lichtenstein's method~\cite{wegmann,birkhoff}, which is also based on the
double-layer (Neumann) kernel, and methods based on the Kerzman-Stein
kernel~\cite{kerzman,murid}. Like our method, these methods are based on an
integral equation whose solution is an easily invertible function of the
boundary correspondence, and they can be used for both interior and exterior
mapping. Similarly, these integral equation methods are also suitable for the
Nyström method with trapezoidal rule. Nevertheless, our method differs from
these because it is the only one which we are aware to make use that the images
of powers of the Riemann map under the double-layer operator are
Faber/Faber-Laurent polynomials.

Other integral equation formulations (e.g.\ due to Berrut~\cite{berrut},
Warschawski~\cite{berrut}, or Banin~\cite{henrici}) produce the derivative of
the boundary correspondence, from which the boundary correspondence may be
recovered by numerical integration. The integral equations of
Gershgorin~\cite{henrici} and
Kantorovich-Krylov~\cite{kantorovich-krylov,henrici} may be used to recover the
boundary correspondence directly, but since the solution is not periodic, this
requires somewhat careful numerical treatment. Perhaps the most well-known
first-kind equation is Symm’s equation, which has been applied to the
computation of both interior and exterior mapping
functions~\cite{symm-interior,symm-exterior}.

For methods for the reverse problem, that is, finding a conformal map from the
interior/exterior of the unit disk onto the interior/exterior of a given domain,
see, for instance~\cite{delillo, delillo-comparison, gutknecht}. For a
comprehensive overview of methods for interior and exterior mapping, including
methods not based on linear boundary integral equations, see
Wegmann~\cite{wegmann}, Gaier~\cite{gaier}, or Henrici~\cite{henrici}. The use
of iterative methods and acceleration techniques also has a long history in the
the solution of systems of equations arising in conformal mapping;
see~\cite{warschawski,trummer,odonnell} for examples.

The remainder of this paper is organized as follows: In
Section~\ref{sec:double-layer}, we recall some facts about harmonic functions
defined on interior or exterior domains, in particular relating to potential
theory and the Cauchy integral, and in Section~\ref{sec:faber}, we discuss the
Faber and Faber-Laurent polynomials. Based on these preliminaries, we introduce
our main technical results regarding the representation of harmonic polynomials
and Laurent polynomials by double-layer potentials in
Section~\ref{sec:dlp-densities}. This allows us to develop a method for interior
and exterior conformal mapping and a high-order discretization method thereof, in
Sections~\ref{sec:mapping-method} and~\ref{sec:numerics}. We close with some numerical experiments on smooth and
non-smooth domains in Section~\ref{sec:num-experiments}.

\section{Preliminaries}

In this paper, we work with a simple, closed, positively oriented curve $\Gamma$,
which we assume to be piecewise smooth. The following conventions will be in
use throughout the paper.

We will refer to the inner component of the curve $\Gamma$ as $\Idom$ and to
the outer component as $\Edom$. Without loss of generality, we will assume $0 \in \Idom$.

We will use $C_R$ to denote the set $\{z : |z| = R\}$.

The interior Riemann map $\Imap$ denotes the complex analytic bijection that
maps $\Idom$ onto $\{ z : |z| < 1 \}$ such that $\Imap(0) = 0$ and $(\Imap)'(0)$
is a positive real number, ensuring uniqueness.

Similarly, we define the exterior Riemann map $\Emap$ as the complex analytic
bijection that maps $\Edom$ onto $\{z : |z| > 1 \}$ for which $\displaystyle
\lim_{z \to \infty} \Emap(z) = \infty $ and $\displaystyle \lim_{z \to \infty}
\pEmap'(z)$ is a positive real number, again ensuring uniqueness.

Carathéodory's theorem~\cite{pommerenke} implies that the interior and exterior
Riemann map continuously extend to the boundary $\Gamma$, establishing a
one-to-one correspondence between $\Gamma$ and the unit circle $C_1$. By the
\emph{boundary correspondences} we will mean the real multi-valued mappings $\Ibcorr$
and $\Ebcorr$, defined on $\Gamma$, such that $\theta^{\pm}(w) = \arg
f^{\pm}(w)$.

\subsection{The Double-Layer Potential}%
\label{sec:double-layer}

The double-layer potential integral operator with density function $\varphi:
\Gamma \to \mathbb{C}$ gives rise to a harmonic function $f$ on the complement
of $\Gamma$
\begin{equation}
\label{eqn:dbl-layer}
f(x) = \mathcal{D} \varphi(x) = - \frac{1}{2 \pi} \int_\Gamma \varphi(y) \,
\dlk{x}{y} \, ds(y), \quad x \in \mathbb{C} \setminus \Gamma.
\end{equation}
Here $\hat{n} \cdot \nabla_y$ denotes the derivative with respect to the
variable $y$ along the outward-facing unit normal $\hat n$ at $y$.
Where the normal is not defined or discontinuous, such as at a corner point, the
kernel has a discontinuity.

Since $x \notin \Gamma$, in a neighborhood of $y$ the value $|y - x|$ is nonzero
and so the logarithm locally possesses a complex analytic branch. The
Cauchy-Riemann equations then imply the relationship
\[
\dlk{x}{y}
= \Re \hat{n} \cdot \nabla_y \log(y - x)
= \Im \hat{\tau} \cdot \nabla_y \log(y - x)
\]
between the normal derivative and the derivative with respect to the
unit tangential vector to the curve, $\hat{\tau}$.
Since
\[
  \hat{\tau} \cdot \nabla_y (\log (y - x))
  =
  \lim_{h \to 0} \frac{1}{h} \left[
    \log \left(y - x + h \hat{\tau} \right)
    - \log(y - x)
    \right]
  = \frac{\hat{\tau}(y)}{y - x},
\]
it follows that the double-layer potential $\mathcal{D}\sigma$ can be written as
\begin{equation}
\label{eqn:dbl-layer-neu-knl}
f(x) = \mathcal{D}\sigma(x) = -\frac{1}{2 \pi} \int_{\Gamma} \sigma(y)
\left( \Im \frac{\hat{\tau}(y)}{y - x} \right)
\, ds(y).
\end{equation}
The kernel appearing in~\eqref{eqn:dbl-layer-neu-knl} is also referred to the
\emph{Neumann kernel}~\cite[Def.~15.9-4]{henrici}.

The equality~\eqref{eqn:dbl-layer-neu-knl} establishes a relationship
between the double-layer potential and the Cauchy integral operator. Since
$\hat{\tau}(y) ds(y) = dy$ and $\Re i\alpha = - \Im \alpha$, we have
\begin{equation}
\label{eqn:dbl-layer-cauchy-knl}
f(x) = -\frac{1}{2 \pi i} \int_\Gamma \sigma(y)
\left( \Re \frac{1}{y - x} \right) dy,
\end{equation}
and thus the kernel of the double-layer potential coincides with the real part
of the Cauchy kernel~\cite[eqn.~(7.37)]{kress}.

\begin{figure}
  \centering
  \begin{tikzpicture}
    \node at (0,0) {%
      \includegraphics{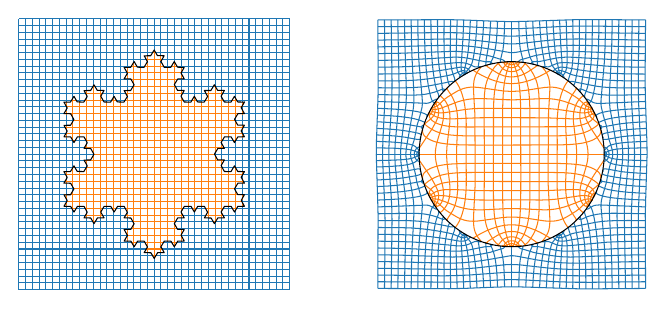}
    };
    \draw [->] (-0.4,0) -- (0.4,0)
      node [pos=0.5,anchor=south] {\color{emapcolor}$\Emap$}
      node [pos=0.5,anchor=north] {\color{imapcolor}$\Imap$};
  \end{tikzpicture}
  \caption{%
    Graphical representation of the Riemann maps for the interior/exterior of a
    Koch snowflake with 192 corners.
  }%
  \label{fig:koch-snowflake}
\end{figure}

\subsection{Faber Polynomials and Faber-Laurent Polynomials}%
\label{sec:faber}

Let $R > 0$ be sufficiently large that the domain $\Idom$ is contained within a
disk of radius $R$ centered at $0$. Then, as $\Emap$ is one-to-one for $|z| > R$, it
follows that $g^+(w):=\Emap(1/w)$ is one-to-one for $0 < |w| < 1/R$.
The pole of $g^+$ at $w = 0$ is simple: Consider the root of $1/g^+$
at the origin. If it were a root of multiplicity $m$ greater than one,
then there would exist an $\epsilon \in \mathbb{C}$ so that the  equation $1/g^+=\epsilon$ has
$m$ simple roots~\cite[Thm.~7.4]{conway_functions_1978}, in contradiction to $g^+$ being
one-to-one. As a result, $g^+$ has the Laurent expansion
\[ g^+(w) = \Emap\Big(\frac1w\Big)=
\frac{\alpha_1}{w} + \alpha_0 + \alpha_{-1} w + \alpha_{-2} w^2 + \cdots,
\quad 0 < |w| < 1/R. \]
This implies that $\Emap$ has the series representation
\begin{equation}
\label{eqn:exterior-map}
\Emap(z) = \alpha_1 z + \alpha_0 +
\frac{\alpha_{-1}}{z} + \frac{\alpha_{-2}}{z^2} + \cdots,
\quad |z| > R.
\end{equation}

The $m$th \emph{Faber polynomial} $p_m(z)$ is defined as the terms of
nonnegative power in the series for $\Emap(z)^m$. It is a polynomial of degree
$m$. Accordingly, ${\Emap(z)}^m$ may be written as
\begin{equation}
  {\Emap(z)}^m = p_m(z) + \hat{p}_m(z)
  \label{eq:emap-pow-m-decomp}
\end{equation}
where $\hat{p}_m(z)$ is a decaying (as $z\to\infty$) function of $z$ defined on $\Edom$.

The $m$th Faber-Laurent polynomial is defined in a similar manner using the
function the $g^-(w):=1/\Imap(w)$ mapping $\Idom$ onto the exterior of the unit
disk. This function, being injective on $\Idom$, has a pole of order 1 at $w =
0$, and so for small enough $r > 0$ admits the Laurent expansion
\[
g^-(w)=
\frac{1}{\Imap(w)} =
\frac{\beta_{-1}}{w} + \beta_0 + \beta_{1} w + \beta_{2} w^2 + \cdots,
\quad 0 < |w| < r. \]
The $m$th \emph{Faber-Laurent polynomial} $q_m$ is defined as the terms of negative power in
the Laurent series for ${g^-(w)}^m$. We have the
representation
\[ {g^-(w)}^m = q_m(z) + \hat{q}_m(z), \]
where $\hat{q}_m(z)$ is a complex analytic function and $q_m(z)$ has a single
order-$m$ pole at $0$.


\section{Representation of Harmonic Polynomials and Laurent Polynomials}%
\label{sec:dlp-densities}

\begin{lemma}%
  \label{lem:happy-lemma}
  Let $m > 0$ be an integer. For all $z \in \Idom$, the $m$th Faber
  polynomial $p_m$ satisfies
  \begin{equation}
  \label{eqn:happy-basis}
  p_m(z) = \frac{1}{\pi} \int_\Gamma {\Emap(y)}^m \, \dlk{z}{y} \, ds(y).
  \end{equation}
\end{lemma}

Effectively, the lemma states that the double-layer potential `filters
out' the decaying part of~\eqref{eq:emap-pow-m-decomp}.

\begin{proof}
  Let $z \in \Idom$.
  Denote the by $I(z)$ the integral
  \[ I(z) = \frac{1}{\pi} \int_\Gamma \Re \left(\Emap(y)^m\right) \, \dlk{z}{y} \, ds(y). \]
  We provide a proof that $I(z) = \Re p_m(z)$, thus handling the case of
  the real part of $p_m$. The argument for the imaginary part of $p_m$ is
  completely analogous to the one for the real part.

  Because $|\Emap(y)| = 1$ for $y \in \Gamma$, we can can write $\Re
  \left(\Emap(y)^m\right)$ as
  \[ \Re \left(\Emap(y)^m\right) = \frac{1}{2} \left( \Emap(y)^m + \frac{1}{\Emap(y)^m} \right). \]
  Using this identity and the identification of the double-layer kernel with the
  real part of the Cauchy kernel (\ref{eqn:dbl-layer-cauchy-knl}), we can
  represent $I(z)$ as
  \begin{equation}
  \label{eqn:cauchy-integral-emap-identification}
  I(z) = \Re \frac{1}{2 \pi i} \int_\Gamma \left( \Emap(y)^m +
  \frac{1}{\Emap(y)^m} \right) \frac{1}{y - z} \, dy.
  \end{equation}

  We proceed by breaking the integral on the right hand side of
  (\ref{eqn:cauchy-integral-emap-identification}) into parts.
  We show
  \begin{equation}
  \label{eqn:cauchy-integral-emap-part1}
  \frac{1}{2 \pi i} \int_\Gamma \frac{\Emap(y)^m}{y - z} \, dy = p_m(z).
  \end{equation}
  We write
  \[ \Emap(y)^m = p_m(y) + \hat{p}_m(y) \]
  where $\hat{p}_m(y)$ goes to $0$ as $|y| \to \infty$. From the Cauchy
  integral formula, we have
  \[ \frac{1}{2 \pi i} \int_\Gamma \frac{p_m(y)}{y - z} \, dy = p_m(z). \]
  To handle $\hat{p}_m$, let $R > 0$ be sufficiently large so that $\Omega$ is
  contained in the interior of the disk with boundary $C_R$. In $\Edom$, the
  function $y \mapsto \hat{p}_m(y)/(y - z)$ is analytic (recall $z \in \Idom$).
  It follows from Cauchy's theorem that
  \[
  \frac{1}{2 \pi i} \int_\Gamma \frac{\hat{p}_m(y)}{y - z} \, dy =
  \frac{1}{2 \pi i} \int_{C_R} \frac{\hat{p}_m(y)}{y - z} \, dy
  \]
  Let $R \to \infty$. Since $\hat{p}_m(y) = O(|y|^{-1})$, the integrand in
  the previous equation is $O(|y|^{-2})$. It follows from a standard integral
  estimate that
  \[
  \lim_{R \to \infty} \left| \frac{1}{2 \pi i} \int_{C_R} \frac{\hat{p}_m(y)}{y
    - z} \, dy \right| = 0,
  \]
  and thus
  \[
  \frac{1}{2 \pi i} \int_\Gamma \frac{\hat{p}_m(y)}{y - z} \, dy = 0.
  \]
  This demonstrates (\ref{eqn:cauchy-integral-emap-part1}) since
  \[
  \frac{1}{2 \pi i} \int_\Gamma \frac{\Emap(y)^m}{y - z} \, dy = \frac{1}{2
    \pi i} \int_\Gamma \frac{p_m(y) + \hat{p}_m(y)}{y - z} \, dy = p_m(z).
  \]

  Next, we show that
  \begin{equation}
  \label{eqn:cauchy-integral-emap-part2}
  \frac{1}{2 \pi i} \int_\Gamma \frac{1}{\Emap(y)^m(y - z)} \, dy = 0.
  \end{equation}
  Again, recall $z \in \Idom$.
  As in the previous paragraph, the integrand is an analytic function of $y$ in
  $\Edom$, so we can choose $R > 0$ sufficiently large so that by Cauchy's
  theorem
  \[
  \frac{1}{2 \pi i} \int_\Gamma \frac{1}{\Emap(y)^m(y - z)} \, dy =
  \frac{1}{2 \pi i} \int_{C_R} \frac{1}{\Emap(y)^m(y - z)} \, dy.
  \]
  Since $\Emap(y)^m = \Theta(|y|^m)$, it follows that the integrand in the
  previous equation is $\Theta(|y|^{-(m+1)})$. By a standard integral estimate,
  we obtain that
  \[
  \lim_{R \to \infty} \left| \frac{1}{2 \pi i} \int_{C_R} \frac{1}{\Emap(y)^m(y
    - z)} \, dy \right| = 0
  \]
  which implies (\ref{eqn:cauchy-integral-emap-part2}).

  By adding together the right hand sides of
  (\ref{eqn:cauchy-integral-emap-part1}) and
  (\ref{eqn:cauchy-integral-emap-part2}) and then taking the real part, we
  obtain that $I(z) = \Re p_m(z)$ via
  (\ref{eqn:cauchy-integral-emap-identification}).
\end{proof}

\begin{lemma}
  \label{lem:happy-lemma-laurent}
  Let $m > 0$ be an integer. For all $z \in \Edom$, the $m$th
  Faber-Laurent polynomial $q_m$ satisfies
  \begin{equation}
  \label{eqn:happy-basis-laurent}
  \overline{q_m(z)} = - \frac{1}{\pi} \int_\Gamma \Imap(y)^m \, \dlk{z}{y} \, ds(y).
  \end{equation}
\end{lemma}

The overline notation $\overline{\,\cdot\,}$ denotes the complex conjugate.

\begin{proof}
  Let $z \in \Edom$ and
  \[ I(z) = -\frac{1}{\pi} \int_\Gamma \Re \left(\Imap(y)^m\right) \, \dlk{z}{y} \, ds(y). \]
  As in the proof of Lemma~\ref{lem:happy-lemma}, we show that $I(z) = \Re
  q_m(z)$, and remark that the imaginary part of $q_m$ can be handled similarly.

  We start by representing $I(z)$ as
  \[
  I(z) = - \Re \frac{1}{2 \pi i}
  \int_\Gamma \left( \Imap(y)^m + \frac{1}{\Imap(y)^m} \right)
  \frac{1}{y - z} \, dy.
  \]

  By Cauchy's theorem, it follows that
  \begin{equation}
  \label{eqn:cauchy-integral-imap-part1}
  - \frac{1}{2 \pi i} \int_\Gamma \frac{\Imap(y)^m}{y - z} \, dy = 0.
  \end{equation}

  We show that
  \begin{equation}
  \label{eqn:cauchy-integral-imap-part2}
  q_m(z) = - \frac{1}{2 \pi i} \int_\Gamma \frac{1}{\Imap(y)^m(y - z)} \, dy.
  \end{equation}
  We proceed by breaking up $1/\Imap(y)^m$ into
  \[ \frac{1}{\Imap(y)^m} = q_m(y) + \hat{q}_m(y) \]
  where $q_m$ is the $m$th Faber-Laurent polynomial, and $\hat{q}_m$ is the
  complex analytic part. First, we handle $q_m$.
  For $y$ close to $0$, we express
  \[
    \frac{1}{y - z} = - \sum_{k=0}^\infty \frac{y^k}{z^{k+1}},
    \qquad
    \text{and}
    \qquad
    q_m(y)=\sum_{k=1}^m \frac{a_k}{y^k}
  \]
  for some $(a_k)$.
  Multiplying both sums, collecting terms, and using the fact that $q_m$ has a
  exactly one pole of order $m$ at $0$, we conclude that
  \[ \Res_{y=0} \left(- q_m(y) \sum_{k=0}^\infty \frac{y^k}{z^{k+1}} \right) = -q_m(z), \]
  which implies by the residue theorem that
  \[ -\frac{1}{2 \pi i} \int_\Gamma \frac{q_m(y)}{y - z} \, dy = q_m(z). \]
  Since $\hat{q}_m$ is an analytic function inside $\Idom$, we have from
  Cauchy's theorem that
  \[ -\frac{1}{2 \pi i} \int_\Gamma \frac{\hat{q}_m(y)}{y - z} \, dy = 0. \]
  This demonstrates~(\ref{eqn:cauchy-integral-imap-part2}).

  The result $I(z) = \Re q_m(z)$ follows by adding together
  (\ref{eqn:cauchy-integral-imap-part1}) and
  (\ref{eqn:cauchy-integral-imap-part2}) and then taking the real part.
\end{proof}

We briefly point out three related results in the literature.  The
basis~\eqref{eqn:happy-basis} can also be derived from~\cite[Lemma 18.2d,
  p.\ 524]{henrici}, although our proof does not rely on this lemma.  Gaier
proves a result similar to the case $m = 1$ of Lemma~\ref{lem:happy-lemma},
in~\cite[p.\ 14, (2.20)]{gaier}, for the case of a horizontal slit. Finally,
in~\cite[(3.12)]{murid} a related integral equation is derived involving the
derivative $\pEmap'$ and the adjoint Neumann kernel.


\section{Integral Equations for Interior and Exterior Conformal Mapping}%
\label{sec:mapping-method}

In this section, we develop a method for recovering the boundary
correspondence assuming that the boundary $\Gamma$ is smooth.

\subsection{Exterior Case}

This section derives an integral equation method for computing the boundary
correspondence $\Ebcorr(z)$ for the exterior map $\Emap$. We solve an integral
equation corresponding to an interior Laplace
Dirichlet problem to obtain a density function $\sigma$, and recover the
boundary correspondence from the density by an application of the Cauchy
integral formula and a normalization.

Recall the power series expansion of the exterior map
\[ \Emap(z) = \alpha_1 z + \alpha_0 + \frac{\alpha_{-1}}{z} +
\frac{\alpha_{-2}}{z^2} + \cdots. \]
From Lemma~\ref{lem:happy-lemma} for $m = 1$, we have that for all $z \in \Idom$
\[
\alpha_1 z + \alpha_0 = \frac{1}{\pi} \int_\Gamma \Emap(y) \, \dlk{z}{y} \, ds(y)
= -2 \mathcal{D} \Emap(z).
\]

Letting $z$ approach a boundary point $\zeta \in \Gamma$ from the interior we obtain, using the
inner jump relation for the double-layer potential~\cite[Thm.~(6.18)]{kress},
the integral equation
\[
\alpha_1 \zeta + \alpha_0 = -2 \mathcal{D} \Emap(\zeta) + \Emap(\zeta),
\quad {\zeta \in \Gamma}.
\]
The parameters $\alpha_1$ and $\alpha_0$ are not assumed to be known a priori.
We use the fact that $\mathcal{D}$ is linear, and that for the constant density
$1(\zeta)$, $\mathcal{D}1(\zeta) = -1/2$ \cite[Ex.~6.17]{kress}. Using these two facts, let us define the
density $\sigma$ as
\[ \sigma(\zeta) = -\frac{1}{\alpha_1}\left( 2\Emap(\zeta) - \alpha_0 \right),
\quad \zeta \in \Gamma. \]
Then we can write the previous integral equation as
\begin{equation}
\label{eqn:exterior-map-ie}
\zeta = \left(\mathcal{D} - \frac{1}{2}\right) \sigma(\zeta), \quad \zeta \in
\Gamma.
\end{equation}
This integral equation is uniquely solvable~\cite[Thm.~6.21]{kress}. From the
density $\sigma$, we can recover $\Emap$ and the boundary correspondence as
follows. As a consequence of~\eqref{eqn:cauchy-integral-emap-part1} established
in Lemma~\ref{lem:happy-lemma}, we have
\[ \frac{1}{2 \pi i} \int_\Gamma \frac{\Emap(y)}{y} \, dy = \alpha_0. \]
Then, from the definition of $\sigma$ and the Cauchy integral formula, we obtain
\[
\frac{1}{2 \pi i} \int_{\Gamma} \frac{\sigma(y)}{y} \, dy
= \frac{1}{2 \pi i} \left[
  -\int_\Gamma \frac{2}{\alpha_1} \cdot \frac{\Emap(y)}{y} \, dy +
  \int_\Gamma \frac{\alpha_0}{\alpha_1} \cdot \frac{1}{y} \, dy \right]
= - \frac{\alpha_0}{\alpha_1}.
\]
Let $\tilde{\sigma}(\zeta)$ denote
\[
\tilde{\sigma}(\zeta) = \sigma(\zeta) +
\frac{1}{2 \pi i} \int_{\Gamma} \frac{\sigma(y)}{y} \, dy
= - \frac{2}{\alpha_1} \Emap(\zeta).
\]
By normalizing $\tilde{\sigma}$ we obtain, for $\zeta \in \Gamma$,
\[
\Emap(\zeta) = -\frac{\tilde{\sigma}(\zeta)}{|\tilde{\sigma}(\zeta)|} \quad
\text{ and } \quad \Ebcorr(\zeta) = \arg \left( -
\frac{\tilde{\sigma}(\zeta)}{|\tilde{\sigma}(\zeta)|} \right).
\]


\subsection{Interior Case}

In this section, we describe a method to recover the interior boundary
correspondence analogous to the exterior one of the previous section. We proceed by describing the solution of an
integral equation corresponding to that of an exterior Laplace Dirichlet
problem for a density function $\sigma$, from which the boundary correspondence
may likewise be recovered by a Cauchy integral and a normalization.

Recall the Laurent series expansion for the inverted interior map
\[
\frac{1}{\Imap(z)} = \frac{\beta_{-1}}{z} + \beta_0 + \beta_1 z + \beta_2 z +
\cdots.
\]
From the case $m = 1$ of Lemma~\ref{lem:happy-lemma-laurent}, we have that for all $z \in
\Edom$,
\[
\overline{\beta_{-1} z^{-1}} =
-\frac{1}{\pi}
\int_\Gamma \Imap(y) \, \dlk{z}{y}
= 2 \mathcal{D}\Imap(z).
\]
Letting $z$ approach a boundary point $\zeta \in \Gamma$ from the exterior, and using the exterior
jump relation for the double-layer potential, we obtain the integral equation
\[
\overline{\beta_{-1} \zeta^{-1}}
= 2 \mathcal{D}\Imap(\zeta) + \Imap(\zeta),
\quad {\zeta \in \Gamma}.
\]
Defining the density $\tilde\sigma(\zeta)$ as
\[
\tilde\sigma(\zeta) = 2(\overline{\beta_{-1}^{-1}}) \Imap(\zeta),
\quad \zeta \in \Gamma,
\]
we can rearrange the above equation to obtain the integral equation
\begin{equation}
\label{eqn:interior-map-ie-nonuniq}
\overline{\zeta^{-1}} =
\left(\mathcal{D} + \frac{1}{2}\right) \tilde\sigma(\zeta),
\quad \zeta \in \Gamma.
\end{equation}
The operator $\left(D + \frac{1}{2}\right)$ has a non-trivial nullspace, which
affects the solvability of this integral equation. Following~\cite{kress}, we remedy
this by defining the operator $\ifunc: C(\Gamma) \to C(\Gamma)$ as
\[ \ifunc\varphi = \int_\Gamma \varphi \, ds. \]
Then the following equation is uniquely solvable~\cite[Thm.~(6.24)]{kress} for a
density $\sigma$:
\begin{equation}
\label{eqn:interior-map-ie-uniq}
\overline{\zeta^{-1}}
= \left(\mathcal{D} + \ifunc{} + \frac{1}{2} \right)
\sigma(\zeta), \quad \zeta \in \Gamma.
\end{equation}

Next, we recall the following facts about the operator $\left( \mathcal{D} +
\frac{1}{2} \right)$. First, the range of $\left( \mathcal{D} + \frac{1}{2}
\right)$ omits the nonzero constant functions. Secondly, the null space of
$\left( \mathcal{D} + \frac{1}{2} \right)$ consists of the
constant functions~\cite[Thm.~(6.21)]{kress}.

If we subtract both sides of~\eqref{eqn:interior-map-ie-uniq} from both sides
of~\eqref{eqn:interior-map-ie-nonuniq}, we obtain that constant function
$\ifunc\sigma$ is in the range of the operator $\left( \mathcal{D} + \frac{1}{2}
\right)$. This implies $\ifunc\sigma = 0$. Thus, we find that
\[
\overline{\zeta^{-1}}
= \left(\mathcal{D} + \frac{1}{2}\right) \tilde\sigma(\zeta)
= \left(\mathcal{D} + \frac{1}{2}\right) \sigma(\zeta), \quad \zeta \in \Gamma.
\]
This implies $\sigma = \tilde\sigma + \delta$ for some $\delta \in \mathbb{C}$.

From the fact that $\Imap(0) = 0$, we know that
\begin{equation*}
  \frac{1}{2 \pi i} \int_{\Gamma} \frac{\sigma(y)}{y} \, dy
  = \frac{1}{2 \pi i} \left[
    2(\overline{\beta_{-1}^{-1}})
    \int_\Gamma
    \frac{\Imap(y)}{y} \, dy +
    \int_\Gamma \frac{\delta}{y} \, dy \right]
  = \delta.
\end{equation*}
Thus we recover $\tilde{\sigma}$ as
\[
\tilde{\sigma}(\zeta)
= \sigma(\zeta) - \frac{1}{2 \pi i} \int_\Gamma \frac{\sigma(y)}{y} \, dy
= 2 (\overline{\beta_{-1}^{-1}}) \Imap(\zeta),
\quad \zeta\in\Gamma.
\]
Recalling $|\Imap(\zeta)|=1$, we normalize to find
\[
\Imap(\zeta) = \frac{\tilde{\sigma}(\zeta)}{|\tilde{\sigma}(\zeta)|} \quad
\text{ and } \quad \Ibcorr(\zeta) = \arg \left(
\frac{\tilde{\sigma}(\zeta)}{|\tilde{\sigma}(\zeta)|} \right),
\quad \zeta\in\Gamma.
\]


\subsection{Summary}%
\label{sec:method-summary}

Algorithm~\ref{alg:conformal-map} captures the operational essence of the
previous two sections. The algorithm is not specific to a particular choice of
discretization, for which a broad range of schemes is applicable.  In the next
section, we provide the details for the Nyström discretization scheme with the
trapezoidal rule for concreteness and for the benefit of our numerical experiments.
This scheme has the advantage of being spectrally accurate,
simple to implement, and amenable to acceleration.

\begin{algorithm}
  \begin{algorithmic}
    \REQUIRE{A smooth Jordan boundary $\Gamma$, with $0$ in the interior.}
    \REQUIRE{A boundary sign $s$: $+1$ for exterior, $-1$ for interior.}
    \ENSURE{Computes the boundary correspondence $\theta$.}
    \STATE{\textsc{Stage 1}}
    \STATE{
      Solve the following integral equation for the density $\sigma$, for all $\zeta \in \Gamma$:
      \[
      \begin{dcases}
        \zeta = \left( \mathcal{D} - \frac{1}{2} \right) \sigma(\zeta) & \text{ if } s = +1 \\
        \overline{\zeta^{-1}} = \left( \mathcal{D} + \ifunc + \frac{1}{2} \right) \sigma(\zeta) & \text{ if } s = -1. \\
      \end{dcases}
      \]}
    \STATE{\textsc{Stage 2}}
    \STATE{Let $\displaystyle \tilde{\sigma}(\zeta) = \sigma(\zeta) + \frac{s}{2 \pi i} \int_\Gamma \frac{\sigma(y)}{y} \, dy$
    ($\zeta\in\Gamma$).}
    \STATE{\textsc{Stage 3}}
    \STATE{Let $\displaystyle \theta(\zeta) =  \arg \left(-s \frac{\tilde{\sigma}(\zeta)}{|\tilde{\sigma}(\zeta)|} \right)$
    ($\zeta\in\Gamma$).}
  \end{algorithmic}
  \caption{Computational method for the obtaining the boundary correspondence}%
  \label{alg:conformal-map}
\end{algorithm}

\section{Numerical Realization of the Methods}%
\label{sec:numerics}

Our main concern in the numerical treatment of Algorithm~\ref{alg:conformal-map}
is the rapid and accurate solution of the integral equations involved.

\subsection{Nyström Method}%
\label{sec:nystrom-method}

We assume a boundary parametrization $\gamma: [0,L] \to \mathbb{C}$ that is $m +
2$ times continuously and periodically differentiable. The operator
$\mathcal{D}$ may by substitution be evaluated in the interval $[0,L]$ using the
parametric Neumann kernel~$\nu$,
\begin{equation}
  \label{eqn:parametric-dlp}
  \mathcal{D} \sigma[\gamma(x)] = \int_0^L \sigma(\gamma(y)) \nu(x, y) \, dy
\end{equation}
which is given by~\cite[p.~394]{henrici}
\[
\nu(x, y) =
-\frac{1}{2 \pi} \Im
\begin{dcases}
  \left. \gamma'(y) \middle/ \left(\gamma(y) - \gamma(x) \right) \right. & x \neq y, \\
  \left. \gamma''(x) \middle/ \left(2 \gamma'(x) \right) \right. & x = y.
\end{dcases}
\]
We consider the discretization of this integral on an $n$ point quadrature rule
with weights ${\{w_j\}}_{j=1}^n$ and nodes $\{y_j\}_{j=1}^n$ on $[0,L]$, which
is given by the functional
\( Q_n g = \sum_{j=1}^n w_j g(y_j). \)
Our specific choice of quadrature rule is the \emph{periodic trapezoidal rule},
whose weights are given by $w_j = L / n$ and the nodes are $y_j = Lj / n$, $j =
1, \ldots, n$.

The Nyström approximation $Q_n[\mathcal{D}]$ to the operator $\mathcal{D}$ uses
pointwise values of the density $\mu = \sigma \circ \gamma$ at the quadrature
nodes as its discrete degrees of freedom,
\[
  Q_n[\mathcal{D}]\mu(x) = \sum_{i=1}^n w_j  \mu(y_j) \nu(x, y_j).
\]
To solve the integral equation
\(
\left(\mathcal{D} - \frac{1}{2}\right) \sigma = f,
\)
we reduce the continuous system to the linear system in $n$ unknowns discretized
at the quadrature points
\begin{equation}
  \label{eqn:nystrom-system}
  Q_n[\mathcal{D}] \mu_n(y_j) - \frac{1}{2} \mu_n(y_j) = f(y_j), \quad j = 1, \ldots,
n.
\end{equation}
Given values of a solution $\mu_n(y_1), \mu_n(y_2), \ldots, \mu_n(y_n)$ to this
system, we may extend $\mu_n$ to a continuous function $\mu_n: [0,L] \to
\mathbb{C}$, by way of the interpolation formula
\begin{equation}
  \label{eqn:nystrom-interpolant}
  \mu_n(x) = 2 \left( Q_n[\mathcal{D}] \mu_n(x) - f(x) \right), \quad x \in [0,L].
\end{equation}
Then, under broadly applicable assumptions on the quadrature rule,
the sequence $Q_n[\mathcal D]$ of operator approximations is invertible for sufficiently
large $n > 0$~\cite[Thm.~12.8]{kress}, has uniformly bounded condition
number~\cite[Thm.~14.3]{kress}, and the sequence $\mu_n$ of approximations to
the density converges uniformly as $n \to \infty$ to the solution of the
continuous system~\cite[Cor.~12.9]{kress}. Furthermore, it can also be
shown~\cite[Cor.~10.14]{kress} the the error in the discrete solution is bounded
from above in the form $\| \mu_n - \mu \|_\infty \leq K \|(Q_n[\mathcal{D}] -
\mathcal{D}) \mu\|$, where $K$ is a constant independent of $n$.

For an $m$ times continuously and periodically differentiable integrand $g$, the
trapezoidal rule admits a spectral error estimate of the form $|Q_n g - \int_0^L
g \, dy| \leq C n^{-m} \| g^{(m)} \|_\infty$, with a constant $C$ independent of
$n$.  It follows that because of the spectral convergence of the periodic
trapezoidal rule we expect spectral convergence in the number of discretization
points for smooth geometries.

\subsection{Fast Iterative Solution of the System}%
\label{sec:accel}

The explicit formation of the dense matrices associated with the
system~\eqref{eqn:nystrom-system} may be avoided by using an iterative method
such as GMRES. Using a Nyström approximation, the number of GMRES iterations for a
fixed accuracy is independent of the number of
unknowns~\cite[Sec.~14.4]{kress}. The iterative application of the operator
$\mathcal{D}$ may be accelerated by considering the discrete operator
$\mathcal{D}$ as the potential due to a set of sources in $\mathbb{R}^2$, and
using the Fast Multipole Method (FMM~\cite{carrier:fmm}). Specifically, we use
the potential
\[
  Q_n[\mathcal{D}] \sigma(y_k) = \frac{1}{2} \omega_k \kappa(y_k) \sigma(y_k) +
  \sum_{\substack{j=1 \\ j \neq k}}^n \omega_j \, \dlk{y_k}{y_j} \, \sigma(y_j), \quad
  {k = 1, \ldots, n}
\]
where $\kappa$ denotes the signed curvature and $\omega_j =
-L/(2 \pi n) \left|\gamma'(y_j)\right|$. On non-pathological particle distributions, the
evaluation phase of the FMM runs in $O(n)$ time and the setup phase takes $O(n
\log n)$ time. It follows that the overall complexity of the solve is $O(n \log
n)$.

\subsection{Evaluation of the Cauchy Integral}%
\label{sec:qbx}

In order to recover the off-boundary values of the Riemann map, one may employ
the Cauchy integral formula for the interior and exterior case
(e.g.~\cite[Eqn.~2.6]{barnett-wu}). At target points $z \in \mathbb{C} \setminus
\Gamma$ sufficiently far from the boundary, quadrature with the trapezoidal rule
is sufficient to achieve high accuracy. However, numerical evaluation of Cauchy
integrals presents challenges close to the boundary $\Gamma$ for standard smooth
quadrature rules such as the trapezoidal rule, leading to the need for an
unacceptably large amount of discretization points~\cite{barnett-wu}. The root
cause of the challenges is the near-singularity of the integrand. A more
efficient strategy for close evaluation is Quadrature by Expansion
(QBX,~\cite{barnett,kloeckner}), a quadrature scheme that exploits the
smoothness of the potential to recover high accuracy near the boundary via
appropriately placed local expansions.   QBX operates by computing
approximate Taylor coefficients of a potential $g$ about centers $c$,
so that an error estimate composed of truncation and quadrature
contributions
\[
\left| g(x) -
\sum_{k=0}^m \frac{\tilde{g}^{(k)}(c)}{k!} (x - c)^k \right| \leq C_1
\|\sigma\|_{C^{m+1}} r^{m+1} + C_2 \left(\frac{h}{4r}\right)^{2q} \| \sigma
\|_{C^{2q}}
\]
can be obtained, dependent on expansion radius $r$, truncation order
$m$, mesh resolution $h$ and quadrature order $q$. We refer
to~\cite{epstein:qbx-error-est} for details.
QBX-based layer potential evaluation may be accelerated by
ways of a fast algorithm~\cite{wala_fast_2018,rachh}, with error
contributions from acceleration very similar to those of conventional
point-based FMMs.
\section{Experimental Results}%
\label{sec:num-experiments}

We implement the method of Section~\ref{sec:mapping-method} using
the numerical approach of Section~\ref{sec:numerics}. In particular, our
discretization is based on the periodic trapezoidal rule with the parametric
Neumann kernel. We employ an FMM-accelerated GMRES solver for the solution of
the integral equations based on FMMLIB~\cite{fmmlib}. The visualizations in
Figures~\ref{fig:koch-snowflake} and~\ref{fig:whale}
were obtained using QBX for evaluation
close to the curve $\Gamma$.

\subsection{Smooth Domains}

\begin{resultstable}
  \caption{Absolute $\ell^\infty$ errors in the boundary
  correspondence on the exterior of an oval of Cassini with shape
  parameter $\alpha$ discretized with $n$ points,
  computed with the methods of Section~\ref{sec:method-summary}
  and~\ref{sec:numerics}.}%
  \label{tab:oval-of-cassini-exterior}
  {\small\begin{tabular}{rSSSSSS}
\toprule
\multicolumn{1}{c}{$n$} & \multicolumn{1}{c}{$\alpha=5$} & \multicolumn{1}{c}{$\alpha=2$} & \multicolumn{1}{c}{$\alpha=1.25$} & \multicolumn{1}{c}{$\alpha=1.11$} & \multicolumn{1}{c}{$\alpha=1.0101$} & \multicolumn{1}{c}{$\alpha=1.001001$}\\
\midrule
$8$ & 2.38e-09 & 2.46e-05 & 4.43e-03 & 2.74e-02 &  & \\
$16$ & 6.96e-16 & 4.66e-09 & 5.44e-05 & 1.25e-03 &  & \\
$32$ &  & 6.28e-16 & 1.15e-08 & 4.23e-06 & 1.57e-02 & \\
$64$ &  &  & 1.05e-15 & 6.79e-11 & 4.15e-04 & 9.82e-02\\
$128$ &  &  &  & 4.00e-16 & 4.73e-07 & 5.90e-03\\
$256$ &  &  &  &  & 8.61e-13 & 7.03e-05\\
$512$ &  &  &  &  & 4.09e-16 & 1.52e-08\\
$1024$ &  &  &  &  &  & 1.11e-15\\
\bottomrule
\end{tabular}
}
\end{resultstable}

\begin{resultstable}
  \caption{Absolute $\ell^\infty$ errors in the boundary
  correspondence on the interior of an epitrochoid with shape
  parameter $\alpha$ discretized with $n$ points,
  computed with the methods of Section~\ref{sec:method-summary}
  and~\ref{sec:numerics}.}%
  \label{tab:epitrochoid-interior}
  {\small\begin{tabular}{rSSSSSS}
\toprule
\multicolumn{1}{c}{$n$} & \multicolumn{1}{c}{$\alpha=0.3$} & \multicolumn{1}{c}{$\alpha=0.4$} & \multicolumn{1}{c}{$\alpha=0.6$} & \multicolumn{1}{c}{$\alpha=0.8$} & \multicolumn{1}{c}{$\alpha=0.9$} & \multicolumn{1}{c}{$\alpha=0.99$}\\
\midrule
$8$ & 2.68e-06 & 2.47e-05 & 1.39e-03 & 3.56e-02 &  & \\
$16$ & 2.55e-12 & 5.38e-10 & 1.80e-06 & 1.24e-03 & 2.41e-02 & \\
$32$ & 3.24e-15 & 2.87e-15 & 2.37e-12 & 3.68e-07 & 1.66e-04 & \\
$64$ &  &  & 1.44e-14 & 5.82e-13 & 5.84e-07 & \\
$128$ &  &  &  & 2.56e-14 & 3.06e-13 & 6.41e-03\\
$256$ &  &  &  &  & 3.63e-14 & 1.80e-04\\
$512$ &  &  &  &  &  & 3.87e-07\\
$1024$ &  &  &  &  &  & 4.66e-12\\
\bottomrule
\end{tabular}
}
\end{resultstable}

We test our method on a number of
smooth test geometries for the interior and exterior case for which the
interior or exterior boundary correspondences are available as analytical expressions. To test the
accuracy of our method, we use the Nyström interpolation formula to evaluate
the compute boundary correspondence at \nevalpts~points on the boundary
equispaced in the parametrization variable, and report the $\ell_\infty$ norm of
the error.

\subsubsection{Oval of Cassini}

\begin{figure}[ht]
  \begin{minipage}[t]{0.5\textwidth}
    \vspace{0pt}
    \centering
    \includegraphics{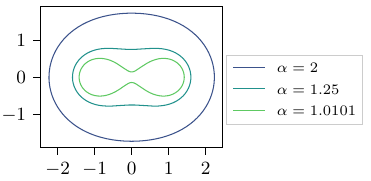}%
    \vspace*{-0.25cm}%
    \caption{\label{fig:cassini-ovals}{Oval of Cassini for various $\alpha$}.}
  \end{minipage}
  \hfill
  \begin{minipage}[t]{0.5\textwidth}
    \vspace{0pt}
    \centering
    \includegraphics{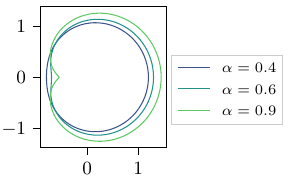}%
    \vspace*{-0.25cm}%
    \caption{\label{fig:epitrochoids}{Epitrochoid for various $\alpha$}.}
  \end{minipage}
\end{figure}

The oval of Cassini curve family is parametrized by $\alpha > 1$. At the
limiting value $\alpha = 1$ the interior is disconnected into two
components. For $\alpha \to \infty$ the domain resembles a disk of radius
$\alpha$. See Figure~\ref{fig:cassini-ovals} for a visualization. The boundary
parametrization and analytical value of the exterior mapping function are given
by:
\begin{align*}
  \gamma_\alpha(t) &= \left(\cos(2t) + \sqrt{a^4 - \sin^2(2t)}\right)^{1/2}
  \exp(it), \quad t \in [0, 2 \pi], \\
  \Emap(\gamma_\alpha(t)) &= (\gamma_\alpha(t)^2 - 1)^{1/2} / \alpha.
\end{align*}
Numerical results demonstrating the accuracy of our method for the
oval of Cassini are shown in Table~\ref{tab:oval-of-cassini-exterior}.

\subsubsection{Epitrochoid}

The epitrochoid family is parametrized by $0 \leq \alpha \leq 1$. At $\alpha =
0$ the boundary is the unit circle, while for $\alpha = 1$ the boundary is a
cardioid. A subset of the tested geometries are visualized in
Figure~\ref{fig:epitrochoids}. The boundary parametrization and the analytical
value of the interior mapping function are given by:
\begin{align*}
  \gamma_\alpha(t) &= \exp(it) + \frac{\alpha}{2} \exp(2it), \quad t \in [0, 2 \pi],
    \\
  \Imap(\gamma_\alpha(t)) &= \exp(it).
\end{align*}
Numerical results demonstrating the accuracy of our method for the epitrochoid are
shown in Table~\ref{tab:epitrochoid-interior}.

\subsubsection{Fourier Whale}

\begin{figure}[ht]
  \begin{minipage}[t]{0.45\textwidth}
    \vspace{0pt}
    \centering
    \includegraphics{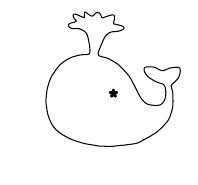}
    \vspace*{-0.5cm}
    \caption{Whale test geometry, analytically
        represented as a Fourier series. The origin is marked.}
    \label{fig:whale-with-origin}
  \end{minipage}
  \hfill
  \begin{minipage}[t]{0.5\textwidth}
      \vspace{0pt}
      \centering
      \resultssetup%
      \captionof{table}{%
        Absolute $\ell^\infty$ self-convergence errors in the boundary correspondence
        for the `Fourier whale' geometry discretized with $n$ points,
        computed with the methods of Section~\ref{sec:method-summary}
        and~\ref{sec:numerics}.}%
      \label{tab:whale}
      {\small\begin{tabular}{rSS}
\toprule
\multicolumn{1}{c}{$n$} & \multicolumn{1}{c}{Interior} & \multicolumn{1}{c}{Exterior}\\
\midrule
$128$ & 7.21e-02 & 1.07e-02\\
$256$ & 2.29e-03 & 4.16e-04\\
$512$ & 1.24e-05 & 3.61e-06\\
$1024$ & 8.10e-10 & 7.90e-11\\
$2048$ & 2.77e-12 & 4.07e-13\\
$4096$ & 2.80e-12 & 4.11e-13\\
\bottomrule
\end{tabular}
}
  \end{minipage}
\end{figure}

\begin{figure}
  \centering
  \begin{tikzpicture}
    \node at (0,0) {%
      \includegraphics{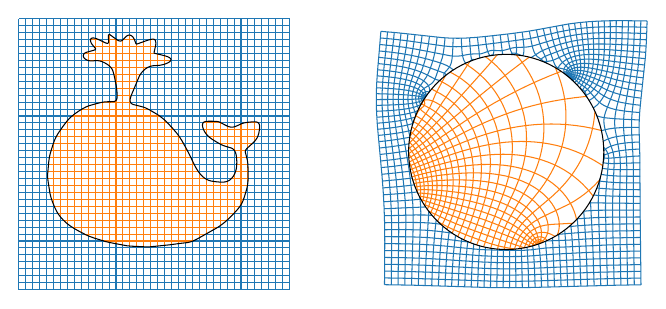}
    };
    \draw [->] (-0.4,0) -- (0.4,0)
      node [pos=0.5,anchor=north] {\color{imapcolor}$\Imap$}
      node [pos=0.5,anchor=south] {\color{emapcolor}$\Emap$};
  \end{tikzpicture}
  \caption{%
    Graphical representation of the Riemann maps for the interior/exterior of the whale
    domain.}%
  \label{fig:whale}
\end{figure}

\begin{figure}[ht]
  \begin{minipage}[t]{0.48\textwidth}
    \vspace{0pt}
    \centering
    \includegraphics{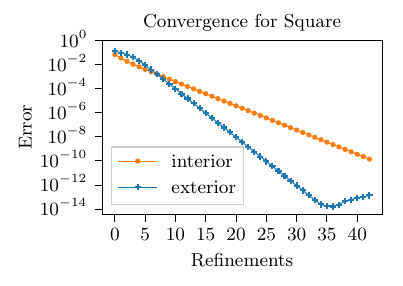}%
    \vspace*{-0.5cm}%
    \caption{Estimated absolute $\ell^\infty$ errors for the unit square with
      increasing refinement. The initial system size was
      $\pointspersquareinitial$ unknowns and each refinement added
      $\pointspersquarerefinement$ unknowns, up to $\pointspersquarefinal$
      unknowns.}%
    \label{fig:square-convergence}
  \end{minipage}
  \hfill
  \begin{minipage}[t]{0.48\textwidth}
    \vspace{0pt}
    \centering
    \includegraphics{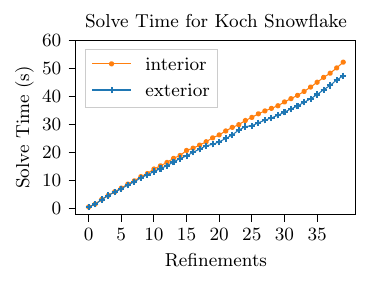}%
    \vspace*{-0.5cm}%
    \caption{Solve time for system associated with the Koch snowflake using an
      FMM-accelerated GMRES solver. The initial system size was
      $\pointspersnowflakeinitial$ unknowns and each refinement added
      $\pointspersnowflakerefinement$ unknowns, up to
      $\pointspersnowflakefinal$ unknowns.
    }%
    \label{fig:koch-snowflake-time}
  \end{minipage}
\end{figure}

We also report the results for a complicated smooth geometry for which the
boundary correspondence is not analytically available. The geometry in
Figure~\ref{fig:whale-with-origin}, with the origin marked, was obtained by
parametrizing the boundary of the image of a spouting whale from the EmojiOne
project~\cite{emojiwhale}. Specifically, we used a parametrization
\(
\gamma(t) = \gamma_1(t) + i\gamma_2(t), t \in [0, 2 \pi],
\)
such that the functions $\gamma_1$ and $\gamma_2$ are given by 53-term Fourier
interpolants of selected boundary points. Since the boundary correspondence is
not analytically available, we test the accuracy using a self-convergence test.
We solve for the boundary correspondence using
$\nwhalesourcepoints$ discretization points.
Using this value as a reference solution, we
estimate the accuracy for a given number of discretization points
by comparing values at $\nevalpts$ points equispaced in
the parameter domain. We report the absolute $\ell^\infty$ error for different
discretization point counts in Table~\ref{tab:whale}. The Riemann maps we
found for this geometry are visualized in Figure~\ref{fig:whale}.

\subsubsection{Discussion: Accuracy on Smooth Geometries}

It is evident from these results that our methods, combined with a
Nyström/trapezoidal scheme, exhibit spectral convergence.
The examples in this section have been used in
previous research to test the accuracy of integral equation methods for
obtaining the boundary correspondence; for instance, for analogous experiments
using the Nyström method applied to the Kerzman-Stein integral equation,
see~\cite{trummer} for the interior case and~\cite{murid} for the exterior
case. We find that these results have a similar level of accuracy to the results
for the Nyström solution of integral equations based on the Kerzman-Stein
kernel. In the next section, we will apply our method to a domain with a corner,
for which the trapezoidal rule is not an ideal quadrature method.

\subsection{Domains with Corners}%
\label{sec:corners}

It is well known that the Nyström method with the trapezoidal rule does not
retain its high order accuracy on boundaries with corners, due to the presence
of singularities in the integral kernel at the corner. Nevertheless, a high
accuracy solution can be recovered with modifications to the scheme. We review
one of these techniques and give numerical results in the case of the square.

The technique we use is based on a simplified version of the quadrature
technique described in~\cite{bremer_nystrom_2011}. We use composite
Gauss-Legendre panels refined dyadically towards the corners of the square. The
last two panels on either side of each corner are omitted. Additionally, each
unknown is multiplied by the square root of the quadrature weight as a way to
improve conditioning (see~\cite{bremer_nystrom_2011} for details).

For the reference solution, we compute the boundary correspondence at 36
equispaced points on the unit square. We use the SC
Toolbox~\cite{driscoll} to compute the Schwarz-Christoffel map from the disk to the square,
and then we invert this map with the provided \texttt{evalinv} subroutine.

The experimental results on the accuracy of our method are given in
Figure~\ref{fig:square-convergence}. The results show the absolute $\ell^\infty$
error in the boundary correspondence versus the number of refinements that were
made recursively to the panels closest to the corner points. The Gauss-Legendre
panels had \squareorder~points per panel. Starting with one panel per side, each
refinement added two panels per side half the width of the previous near-corner panels.

This scheme is able to recover up to 13--14 digits of accuracy in the exterior
case and 10 digits in the interior case.  The difference in convergence speeds
between the interior and the exterior of the square is the subject of future
investigation. We have chosen this scheme for its simplicity, though more
advanced schemes can both improve the accuracy and reduce the number of unknowns
required.

\subsection{Scaling}

To study the scaling of our method, we time our implementation of the solve
phase of Algorithm~\ref{alg:conformal-map}, which is the dominant contribution
to the cost of the algorithm. The test geometry is the fourth iteration of a
Koch snowflake curve, which has 192 corners (see Figure~\ref{fig:koch-snowflake}
for a visualization) and the system contains up to
$\pointspersnowflakefinal$ unknowns. We use the quadrature scheme described in
Section~\ref{sec:corners} and measure the wall time of the algorithm with
increasing refinements. The timing results, obtained on a single core
of a dual-socket 2.2~GHz
Intel~Xeon~E5\nobreakdash-2650~v4 processor, are presented in
Figure~\ref{fig:koch-snowflake-time}. After 39 refinements, the mapping obtains
approximately 7--8 digits of accuracy according to direct comparison with
results from the SC Toolbox. As expected, the timing data demonstrates the solve
phase of the algorithm scales close to linearly with the number of unknowns.

\section{Conclusions}

This paper makes two contributions.

First, we characterize the density functions $\sigma$ that give rise to
harmonic polynomials represented as double-layer potentials $\mathcal{D} \sigma$
on the interior of a piecewise smooth Jordan domain, and their counterparts that
give rise to Laurent polynomials on the exterior of the domain. We show how
these density functions relate to the Riemann maps associated with the
domain. In addition to the described application to conformal mapping, this work
may be of mathematical interest for those studying the behavior of the
double-layer potential and numerical methods for it.

Second, we derive an integral equation whose solution allows us to recover the
boundary correspondence for the exterior or interior mapping function. From a
practical standpoint, our equation is second-kind, uniquely solvable and has a
continuous kernel, which leads to a robust and simple discretization with the
Nyström method.
We further demonstrate the effective acceleration of the method,
avoiding quadratic complexity in all parts of the method.
A major advantage of the double-layer potential is the ready
availability of existing fast solvers, such as those
in~\cite{bremer,bremer-gillman,ho,rachh,smigaj}. Our experiments demonstrate
that the method achieves spectral accuracy on smooth domains, with results of
comparable accuracy to those based on the solution of integral equations using
the Kerzman-Stein kernel. Finally, we demonstrate that the method can be made accurate in the
presence of corners.


\section*{Acknowledgments}

The authors' research was supported by the National Science Foundation under
grants DMS-1418961 and DMS-1654756.
Part of the work was performed while the authors were participating in
the HKUST-ICERM workshop `Integral Equation Methods, Fast
Algorithms and Their Applications to Fluid Dynamics and Materials
Science' held in 2017. The authors would also like to thank the anonymous
reviewers whose comments have helped to improve the manuscript.
\bibliography{main}{}
\bibliographystyle{siam}

\end{document}
